# Split Equilibrium Problems for Related Games

# and Applications to Economic Theory

Jinlu Li

Department of Mathematics
Shawnee State University
Portsmouth, Ohio 45662
U. S. A.

**Abstract**: In this paper we introduce the concept of split Nash equilibrium problems associated with two related noncooperative strategic games. Then we apply the Fan-KKM theorem to prove the existence of solutions to split Nash equilibrium problems of related noncooperative strategic games, in which the strategy sets of the players are nonempty closed and convex subsets in Banach spaces. As application of this existence to economics, an example is provided to study the existence of split Nash equilibrium of utilities of two related economies. As applications, we study the existence of split Nash equilibrium in the repeated extended Bertrant duopoly model of price competition

**Keywords**:  Fan-KKM theorem; split Nash equilibrium problem; convexity direction preserved mapping; extended Bertrant duopoly model of price competition

**2010 Mathematics Subject Classification:** 49J40; 49J52; 91A10; 91A25; 91A80

## 1. Introduction

Nash equilibrium theory has become an important branch in game theory (see [1], [12], [19]). Together with Pareto equilibrium theory, Nash equilibrium theory has been widely applied to many different fields, particularly in economics (see [9], [11], [15], [16], [19]).

In pure and applied mathematics, Nash equilibrium problems have attracted numerous authors to study, leading to its generalization and extension to broad fields with many applications (see [5], [11], [17]).

Recently, along with the development of split variational inequality theory, many authors have been rapidly developing the split equilibrium problem theory (see [2-5], [8-9], [14], [18]). These authors studied mathematical optimization problems which are considered to be equilibrium problems. To my knowledge, there have been no authors studying split Nash equilibrium problems for noncooperative strategic games.

Since the source of Nash equilibrium problems is noncooperative strategic games, we consider split Nash equilibrium problems for two related noncooperative strategic games (That is defined in section 3). In this paper, we introduce the concept of split Nash equilibrium problems associated with two noncooperative strategic games with the strategy sets for the players to be nonempty closed and convex subsets of Banach spaces. We utilize the Fan-KKM theorem to prove the existence of solutions to some split Nash equilibrium problems of noncooperative strategic games.

The results of the split Nash equilibrium problems of noncooperative strategic games can be directly applied to the existence of split Nash equilibrium for repeating noncooperative strategic games in which the strategies in the repeated games are linearly transformed from the previous games.

This paper is organized as follows: for reference, in section 2, we recall some concepts of partially ordered topological vector spaces, the concept of KKM mapping and the Fan-KKM theorem; in section 3, we introduce the concept of split Nash equilibrium problems and the concept of convexity direction preserved property for linear operators. Then by applying the Fan-KKM theorem, we prove an existence theorem of solutions to some split Nash equilibrium problems of related noncooperative strategic games. To my knowledge, this is the first time an existence theorem for split Nash equilibrium problems has been proven. As application of this existence to economics, an example is provided to study the existence of split Nash equilibrium of utilities of two related economies. Then, we provide an existence theorem of split Nash equilibrium problems for repeated noncooperative strategic games, that can be considered as a special case of split Nash equilibrium problems of related noncooperative strategic games. Finally, we generalize the Bertrant duopoly model to the extended Bertrant model of price competition and consider the Markov split Nash equilibrium problems in the repeated extended Bertrant duopoly model.

## 2. Preliminaries

### 2.1. Partially ordered topological vector spaces

Suppose that a topological vector space $U$ is equipped with a partial order $\succcurlyeq$ which has the following properties:

1. $w \succcurlyeq u$ implies $w + z \succcurlyeq u + z$, for all $u, w, z \in U$;
2. $w \succcurlyeq u$ implies $\alpha w \succcurlyeq \alpha u$, for $u, w \in U$ and $\alpha \geq 0$;
3. For any $u \in U$, the following sets are closed

$$[u) = \{w \in U: w \succcurlyeq u\} \text{ and } (u] = \{w \in U: u \succcurlyeq w\}.$$

Then, in this paper, $U$ is called a partially ordered topological vector space, or is simply called an ordered topological vector space, which is written as $(U, \succcurlyeq)$. The above Condition 3 of partially ordered topological vector spaces is the definition that the topology of $U$ is natural with respect to the partial order $\succcurlyeq$ on $U$.

### 2.2. Fan-KKM Theorem

Fan-KKM theorem is applied in the proof of the main theorem in this section. For more details regarding to KKM mappings and Fan-KKM theorem, please refer to Fan [6] and Park [13]. For easy reference, we briefly recall the definition of KKM mappings and Fan-KKM theorem below.

Let $K$ be a nonempty subset of a linear space $B$. A set-valued mapping $T: K \to 2^B \backslash \{\emptyset\}$ is said to be a KKM mapping if for any finite subset $\{y_1, y_2, \ldots, y_n\}$ of $K$, we have

$$\text{co}\{y_1, y_2, \ldots, y_n\} \subseteq \cup_{1 \leq i \leq n} T(y_i),$$

where co$\{y_1, y_2, \ldots, y_n\}$ denotes the convex hull of $\{y_1, y_2, \ldots, y_n\}$.

**Fan-KKM Theorem**. *Let $K$ be a nonempty closed convex subset of a Hausdorff topological vector space $B$ and let $T: K \to 2^B \backslash \{\emptyset\}$ be a KKM mapping with closed values. If there exists a point $y_0 \in K$ such that $T(y_0)$ is a compact subset, then*

$$\cap_{y \in K} T(y) \neq \emptyset.$$

## 3. Split Nash equilibrium problems in related noncooperative strategic games

### 3.1. Noncooperative strategic games

We first recall some concepts of noncooperative strategic games. For the details see [12] and [19]. Let $n$ be a positive integer greater than 1. An $n$-person noncooperative strategic game consists of the following elements:

1. the set of $n$ players denoted by $N$ with $|N| = n$;
2. the set of profiles $S_N = \Pi_{i \in N} S_i$, where $S_i$ is the pure strategy set for player $i \in N$;
3. the utility (payoff) function $f_i$ for $i \in N$, which is a real function defined on $S_N$.

This game is denoted by $G(N, S_N, f)$, where $f = \Pi_{i \in N} f_i$, that is a mapping from $S_N$ to $R^n$. As usual, for every $i \in N$, we often denote a profile of pure strategies for player $i$'s opponents by

$$x_{-i} = (x_1, x_2, \ldots, x_{i-1}, x_{i+1}, x_n).$$

We similarly denote the set of profiles of pure strategies for player $i$'s opponents by

$$S_{-i} = \Pi_{k \in N \backslash \{i\}} S_k.$$

Throughout this section, unless otherwise is stated, we assume that, in the products $\Pi_{i \in N} S_i$, $\Pi_{i \in N} f_i$ and $\Pi_{k \in N \backslash \{i\}} S_k$, the players appear in the same sequential orders. Then we write $x \in S_N$ as

$$x = (x_i, x_{-i}) \text{ with } x_{-i} \in S_{-i}, \text{ for } i \in N. \tag{1}$$

Moreover, for every $x_{-i} \in S_{-i}$, we denote

$$f_i(S_i, x_{-i}) = \{f_i(z_i, x_{-i}): z_i \in S_i\}.$$

From $f = \Pi_{i \in N} f_i$ in game $G(N, S_N, f)$, for any $x \in S$, we have

$$f(x) = \Pi_{i \in N} f_i(x).$$

For a given $f = \Pi_{i \in N} f_i$, we conveniently define a mapping $F: S_N \times S_N \to R^n$ by

$$F(z, x) = \Pi_{i \in N} f_i(z_i, x_{-i}), \text{ for any } x, z \in S_N.$$

If $x = z$, then

$$F(z, x) = F(x, x) = f(x).$$

Let $\geq^n$ be the component-wise ordering on $R^n$. It is clear to see that $\geq^n$ is a partial order on $R^n$. It follows that, for $x, z \in S_N$,

$$F(z, x) \leq^n f(x) \text{ if and only if } f_i(z_i, x_{-i}) \leq f_i(x_i, x_{-i}), \text{ for all } i \in N.$$

In an $n$-person non-cooperative strategic game $G(N, S_N, f)$, a profile $\hat{x} \in S_N$ is called a Nash equilibrium of this game if the following inequality holds:

$$F(z, \hat{x}) \leq^n f(\hat{x}), \text{ for all } z \in S_N.$$

It is equivalent to

$$f_i(z_i, \hat{x}_{-i}) \leq f_i(\hat{x}_i, \hat{x}_{-i}), \text{ for every } i \in N \text{ and for every } z_i \in S_i.$$

In an $n$-person non-cooperative strategic game $G(N, S_N, f)$, for every $i \in N$, we define a set-valued mapping $\Gamma_i: S_i \to 2^{S_N}$ by

$$\Gamma_i(x_i) = \{(z_i, z_{-i}) \in S_N: f_i(x_i, z_{-i}) \leq f_i(z_i, z_{-i})\}, \text{ for all } x_i \in S_i.$$

Then a set-valued mapping $\Gamma_N: S_N \to 2^{S_N}$ is defined by

$$\Gamma_N(x) = \{(z_i, z_{-i}) \in S_N: f_i(x_i, z_{-i}) \leq f_i(z_i, z_{-i}), \text{ for all } i \in N\}, \text{ for } x \in S_N.$$

That is, for every $x \in S_N$,

$$\Gamma_N(x) = \cap_{i \in N} \Gamma_i(x_i) = \cap_{i \in N} \{(z_i, z_{-i}) \in S_N: f_i(x_i, z_{-i}) \leq f_i(z_i, z_{-i})\}.$$

Using the notations for $F(z, x)$ and $F(x, x)$, it can be rewritten as

$$\Gamma_N(x) = \{z \in S_N: F(x, z) \leq^n f(z)\}, \text{ for } x \in S_N.$$

It can be shown that the standard topology on $R^n$ is natural with respect to the partial ordering $\leq^n$. That is, for every $u \in R^n$, the following subsets of $R^n$

$$\{v \in R^n: v \leq^n u\} \text{ and } \{v \in R^n: u \leq^n v\}$$

are all closed. The following lemma immediately follows from the above facts.

**Lemma 3. 1**. *Let* $G = (N, S_N, f)$ *be an n-person non-cooperative strategic game. Suppose that, for every* $i \in N$, $S_i$ *is a nonempty closed and convex subset of a Banach space* $B_i$ *and* $f_i$ *is a continuous real valued function defined on* $S_N$. *Then, for any* $x \in S_N$, $\Gamma_N(x)$ *is a closed subset of* $S_N$ *with respect to the product topology on* $S_N$.

### 3. 2. Split Nash equilibrium Problems for two related noncooperative strategic games

Let $G(M, S_M, g)$ be an $m$-person non-cooperative strategic game. Similar to the game $(N, S_N, f)$ considered in the previous sub-section, in the game $G(M, S_M, g)$, the player set is $M$ with $|M| = m$, for some $m > 1$. For every player $j \in M$, $S_j$ is the strategy set and $g_j$ is his utility (payoff) function which is a real function defined on $S_M$. We also write the set of profiles $S_M = \Pi_{j \in M} S_j$. The product utilities $g = \Pi_{j \in M} g_j$ is a mapping from $S_M$ to $R^m$.

We can similarly define a mapping $G: S_M \times S_M \to R^m$ by

$$G(w, y) = \Pi_{j \in M} g_j(w_j, y_{-j}), \text{ for } y, w \in S_M.$$

Let $G(N, S_N, f)$ and $G(M, S_M, g)$ be $n$-person and $m$-person non-cooperative strategic games, respectively. Suppose that, for every $i \in N$, $S_i$ is a nonempty closed and convex subset of a Banach space $B_i$ and, for every $j \in M$, $S_j$ is a nonempty closed and convex subset of a Banach space $B_j$. We denote the product Banach spaces by

$$B_N = \Pi_{i \in N} B_i \text{ and } B_M = \Pi_{j \in M} B_j.$$

Throughout this section, unless otherwise stated, we assume that the strategy set of every player in a non-cooperative strategic game is a nonempty closed and convex subset of a Banach space.

Let $A \in \mathcal{L}(B_N, B_M)$, that is, $A: B_N \to B_M$ is a linear continuous operator. Suppose that $A$ has the component-wise linear continuous mappings $A_j \in \mathcal{L}(B_N, B_j)$, for $j \in M$, such that

$$A(x) = \Pi_{j \in M} A_j(x) = \Pi_{j \in M} (Ax)_j \in B_M, \text{ for } x \in B_N.$$

Then, similar to the notation in (1), for $j \in M$, we can simply write

$$A(x) = (A_j(x), A_{-j}(x)) = ((Ax)_j, (Ax)_{-j}) = (Ax_j, Ax_{-j}), \tag{2}$$

where $A_{-j}(x) = \Pi_{k \in M, k \neq j} A_k(x) \in B_{-j}$ and $B_{-j} = \Pi_{k \in M, k \neq j} B_k$.

Let $G(N, S_N, f)$ and $G(M, S_M, g)$ be $n$-person and $m$-person non-cooperative strategic games, respectively. If there is a linear continuous operator $A \in \mathcal{L}(B_N, B_M)$ satisfying the following condition (3), then $G(N, S_N, f)$ and $G(M, S_M, g)$ are said to be related.

$$AS_N \subseteq S_M. \tag{3}$$

**Definition 3.2.** Let G($N$, $S_N$, $f$) and G($M$, $S_M$, $g$) be n-person and m-person related non-cooperative strategic games, respectively, with a linear continuous operator A: $B_N \to B_M$. The split Nash equilibrium problem associated with the games G($N$, $S_N$, $f$), G($M$, $S_M$, $g$), and the operator A, denoted by SNE(G($N$, $S_N$, $f$), G($M$, $S_M$, $g$), A), is formalized as:

$$\text{find a Nash equilibrium } \hat{x} \text{ of the game G}(N, S_N, f), \tag{4}$$

such that $\quad A\hat{x}$ is a Nash equilibrium of the game G($M$, $S_M$, $g$). $\tag{5}$

That is, to find a profile $\hat{x} \in S_N$ satisfying

$$f_i(z_i, \hat{x}_{-i}) \leq f_i(\hat{x}_i, \hat{x}_{-i}), \text{ for every } i \in N \text{ and for every } z_i \in S_i, \tag{4'}$$

such that the profile $A\hat{x} \in S_M$ solves

$$g_j(y_j, A\hat{x}_{-j}) \leq g_j(A\hat{x}_j, A\hat{x}_{-j}), \text{ for every } j \in M \text{ and for every } y_j \in S_j. \tag{5'}$$

Recall that, from the notation introduced in (2), for every $j \in M$, we have

$$A\hat{x}_j = (A\hat{x})_j = A_j\hat{x} \text{ and } A\hat{x}_{-j} = (A\hat{x})_{-j} = A_{-j}\hat{x}.$$

It can be equivalently written as: to find a profile $\hat{x} \in S_N$ satisfying

$$F(z, \hat{x}) \leq^n f(\hat{x}), \text{ for all } z \in S_N. \tag{4''}$$

such that the profile $A\hat{x} \in S_M$ solves

$$G(y, A\hat{x}) \leq g_j(A\hat{x}), \text{ for every } y \in S_M. \tag{5''}$$

Such a profile $\hat{x}$ in $S_N$ is called a split Nash equilibrium of this split Nash equilibrium problem SNE(G($N$, $S_N$, $f$), G($M$, $S_M$, $g$), A). The set of all solutions to this split Nash equilibriums in this problem is denoted by

$$\mathcal{S}(\text{G}(N, S_N, f), \text{G}(M, S_M, g), A)).$$

When looking at the equilibrium problems (4) and (5) separately, the problem (4) is the classical Nash equilibriums problem of non-cooperative strategic games. When considering a special case of SNE(G($N$, $S_N$, $f$), G($M$, $S_M$, $g$), A) such as: $N = M$, $S_N = S_M$, $f = g$ (G($N$, $S_N$, $f$) = G($M$, $S_M$, $g$)) and $A = I$ that is the identity mapping on $S_N$, the split Nash equilibrium problem SNE(G($N$, $S_N$, $f$), G($M$, $S_M$, $g$), A) reduces to the classical Nash equilibrium problem for the game G($N$, $S_N$, $f$). In this view, split Nash equilibrium problems can be considered as the natural extensions of the classical Nash equilibrium problems.

### 3.3. Convexity direction preserved linear operators

Let $G(N, S_N, f)$ and $G(M, S_M, g)$ be $n$-person and $m$-person non-cooperative strategic games with concave utility functions, respectively. A linear continuous operator $A: B_N \to B_M$ is said to be *convexity direction preserved* with respect to $f$ and $g$ on $S_N$ and $S_M$ if, for any given points $(u, Au)$, $(v, Av) \in S_N \times S_M$ and for their arbitrary convex combination $w = \lambda u + (1-\lambda)v$, where $0 \leq \lambda \leq 1$, we have either that

$$F(u, w) \leq^n f(w) \text{ and } G(Au, Aw) \leq^m g(Aw) \text{ both hold}$$

or that

$$F(v, w) \leq^n f(w) \text{ and } G(Av, Aw) \leq^m g(Aw) \text{ both hold.}$$

When understanding the concept of convexity direction preserved property, we can see that, for any given points $u, v \in S_N$ and their convex combination of $w = \lambda u + (1-\lambda)v$, we must have

$$\text{either } F(u, w) \leq^n f(w), \text{ or } F(v, w) \leq^n f(w). \tag{6}$$

To see this, we assume by contradiction that

$$F(u, w) >^n f(w) \text{ and } F(v, w) >^n f(w) \text{ both hold.}$$

Then from the concavity of $f_i$'s, we have

$$\begin{aligned} f(w) &= F(w, w) \\ &= F(\lambda u + (1-\lambda)v, \lambda w + (1-\lambda)w) \\ &\geq^n \lambda F(u, w) + (1-\lambda)F(v, w) \\ &>^n \lambda f(w) + (1-\lambda) f(w) \\ &= f(w). \end{aligned}$$

This is a contradiction. Hence (6) must hold. Since $A \in \mathcal{L}(B_N, B_M)$, it follows that

$$Aw = \lambda Au + (1-\lambda)Av.$$

Then we can similarly show that

$$\text{either } G(Au, Aw) \leq^m g(Aw) \text{ or } G(Av, Aw) \leq^m g(Aw). \tag{7}$$

The convexity direction preserved property of $A$ insures that one of the two inequalities in (6) and the corresponding inequality in (7) must simultaneously hold for the point $u$ or the point $v$.

### 3.4. An existence theorem

Now we state and prove the main theorem in this section.

**Theorem 3.3**. *Let* $G(N, S_N, f)$ *and* $G(M, S_M, g)$ *be n-person and m-person related non-cooperative strategic games, respectively, with a linear continuous operator* $A: B_N \to B_M$. *Suppose that* $f_i$'s, $g_j$'s, *and* $A$, *satisfy the following conditions*:

(a1). $f_i$'s *and* $g_j$'s *are all continuous and concave*;
(a2). $AS_N = S_M$;
(a3). $A$ *is convexity vector direction preserved with respect to f and g on* $S_N$ *and* $S_M$.

*If there is a point* $(t, At) \in S_N \times S_M$ *such that*

$$\{(z, Az) \in S_N \times S_M: F(t, z) \leq^n f(z) \text{ and } G(At, Az) \leq^m g(Az)\} \text{ is compact,}$$

*then* SNE($G(N, S_N, f)$, $G(M, S_M, g)$, $A$) *has a solution*.

**Proof**. For this operator $A \in \mathcal{L}(B_N, B_M)$ satisfying (a2) in this theorem, denote $K$ its graph. $K$ is a subset of $S_N \times S_M \subseteq B_N \times B_M$ and is defined by

$$K = \{(x, Ax): x \in S_N\} \subseteq S_N \times S_M.$$

Since $B_N$ and $B_M$ are Banach spaces with the product topologies, by the closed graph theorem in Banach spaces, the graph of the bounded linear operator $A$ from $S_N$ to $S_M$ is a closed subset of $S_N \times S_M$ with respect to the product topology. Then the convexity of $S_N$ and $S_M$ and the linearity of $A$ imply that $K$ is a nonempty closed and convex subset of $S_N \times S_M \subseteq B_N \times B_M$.

On this underlying space $K$, we define a mapping $T: K \to 2^K \setminus \{\emptyset\}$ by

$$T(x, Ax) = \{(z, Az) \in K: F(x, z) \leq^n f(z) \text{ and } G(Ax, Az) \leq^m g(Az)\}, \text{ for } (x, Ax) \in K.$$

From the continuity of the mappings $A$, $f$ and $g$, and (7), we have $(x, Ax) \in T(x, Ax)$ which yields that, for every $(x, Ax) \in K$, $T(x, Ax)$ is a nonempty closed subset of $K \subseteq S_N \times S_M \subseteq B_N \times B_M$.

Next we show that $T$ is a KKM mapping. For any points $(u, Au), (v, Av) \in K$, we arbitrarily take a convex combination of $w = \lambda u + (1-\lambda)v$, where $0 < \lambda < 1$. Since $A \in \mathcal{L}(B_N, B_M)$, it follows that

$$Aw = \lambda Au + (1-\lambda)Av.$$

Then

$$(w, Aw) = \lambda(u, Au) + (1-\lambda)(v, Av)$$

From condition (a2), $A$ has the convexity direction preserved property with respect to $f$ and $g$ on $S_N$ and $S_M$. Then, we have either

$$F(u, w) \leq^n f(w) \text{ and } G(Au, Aw) \leq^m g(Aw),$$

or

$$F(v, w) \leq^n f(w) \text{ and } G(Av, Aw) \leq^m g(Aw).$$

That is,

$$\text{either } (w, Aw) \in T(u, Au), \text{ or } (w, Aw) \in T(v, Av).$$

It implies that

$$(w, Aw) \in T(u, Au) \cup T(v, Av).$$

Similarly, we can extend this to finite convex combinations. Hence, $T$ is a KKM mapping. For the given point $(t, At) \in K$ from the assumptions of this theorem, we have

$$T(t, At) = \{(z, Az) \in K: F(t, z) \leq^n f(z) \text{ and } G(At, Az) \leq^m g(Az)\}.$$

It is a nonempty compact subset of $K \subseteq S_N \times S_M$.

By applying the Fan-KKM theorem, we obtain that $\bigcap_{(x, Ax) \in K} T(x, Ax) \neq \emptyset$. Then, by taking any $(x_*, Ax_*) \in \bigcap_{(x, Ax) \in K} T(x, Ax)$, we have

$$F(x, x^*) \leq^n f(x^*), \text{ for all } x \in S_N, \tag{8}$$

and

$$G(Ax, Ax^*) \leq^m g(Ax^*), \text{ for all } x \in S_N. \tag{9}$$

Let $y_* = Ax_*$. Notice that, for the fixed point $x_* \in S_N$ and $Ax_* \in S_M$, the inclusion properties (8) and (9) can be considered separately. When we look at (9) for the fixed $x_* \in S_N$, from condition (a2) of this theorem, (9) is equivalent to

$$G(y, Ax^*) \leq^m g(Ax^*), \text{ for all } y \in S_M. \tag{10}$$

Combining (8), (10), and $y_* = Ax_*$, we have $(x_*, y_*) = (x_*, Ax_*)$ is a solution of the split Nash equilibrium problem SNE(G($N$, $S_N$, $f$), G($M$, $S_M$, $g$), $A$). □

A consequence of Theorem 3.3:

**Corollary 3.4.** *Let* G($N$, $S_N$, $f$) *and* G($M$, $S_M$, $g$) *be n-person and m-person non-cooperative strategic games, respectively, with the strategy sets all nonempty compact and convex subsets of Banach spaces. Let $A: B_N \to B_M$ be a linear continuous operator. Suppose that $f_i$'s, $g_j$'s and $A$ satisfy the following conditions:*

(a1). *$f_i$'s and $g_j$'s all are continuous and concave;*
(a2). *$AS_N = S_M$;*

(a3). *A is convexity vector direction preserved with respect to f and g on $S_N$ and $S_M$.*

*Then* SNE(G(N, $S_N$, f), G(M, $S_M$, g), A) *has a solution.*

**Remarks 3.5**. The results in Theorem 3.3 and Corollary 3.4 do not mean that there does not exist a split Nash equilibrium if the conditions of these existence results do not hold. It only means that it cannot be assured that there is one.

## 4. An example

In this section we consider the Nash equilibrium problems for the utilities of industries in two related economic forests. In view of game theory, it provides an example of split Nash equilibrium for noncooperative strategic games.

**Example 4.1**. Suppose that there are two related economies, $E_1$ and $E_2$. In economy $E_1$, there are three industries, *a*, *b* and *c*, which produce certain types of goods provided to the market that are defined by variables *x*, *y*, and *z*, respectively. Suppose that the quantities are nonnegative and the utilities (payoffs) are dependable. They are respectively defined by

$$f_a(x, y, z) = xyz - 4x^2; \tag{11}$$

$$f_b(x, y, z) = x^2 yz - \frac{1}{8} y^4; \tag{12}$$

$$f_c(x, y, z) = x^{\frac{1}{2}} y z^{\frac{1}{2}} - \frac{1}{2} z, \quad \text{for } 0 \leq x, y, z < \infty. \tag{13}$$

This economy formalizes a noncooperative strategic game, which has player set $N = \{a, b, c\}$ with $|N| = 3$. For each player, the strategy set is $R_+ = [0, \infty)$. Hence the set of profiles $S_N = R_+^3$. The utility functions $f_a$, $f_b$, and $f_c$, of the players are defined in (11–13) with $f = f_a \times f_b \times f_c$. This game, which represents the economy $E_1$, is denoted by G(N, $S_N$, f).

In economy $E_2$, there are two industries, *d* and *e*, which also produce certain types of goods that are defined by variables *s* and *t*, respectively. Suppose that the quantities are also nonnegative and the utilities are also dependable and respectively defined by

$$g_d(s, t) = \frac{1}{2} st - \frac{1}{3} s^2; \tag{14}$$

$$g_e(s, t) = 48 s^{\frac{1}{2}} t - \frac{1}{48} t^4, \quad \text{for } 0 \leq s, t < \infty. \tag{15}$$

This economy formalizes a noncooperative strategic game, which has player set $M = \{d, e\}$ with $|M| = 2$. For each player, the strategy set is $R_+ = [0, \infty)$. Hence the set of profiles $S_M = R_+^2$. The utility functions $g_d$ and $g_e$ of the players are defined in (14–15) with $g = g_d \times g_e$. This game, which represents the economy $E_2$, is denoted by $G(M, S_M, g)$.

Suppose that the two economies are related by their outcomes and the market demands in $E_2$ are higher than that in $E_1$. As in Section 3, the related relationship between the two economies $E_1$ and $E_2$ is represented by the following operator $A: R_+^3 \to R_+^2$,

$$(s, t)^T = \begin{pmatrix} s \\ t \end{pmatrix} = A \begin{pmatrix} x \\ y \\ z \end{pmatrix},$$

where $A$ is a positive linear and continuous operator from $R_+^3$ to $R_+^2$, which is defined by a 2×3 positive matrix below

$$A = \begin{pmatrix} 1 & 2 & 1 \\ 2 & 1 & 2 \end{pmatrix} = \begin{pmatrix} A_d \\ A_e \end{pmatrix},$$

with $A_d = (1\ 2\ 1)$ and $A_e = (2\ 1\ 2)$ satisfying $A_d, A_e \in \mathcal{L}(R_+^3, R_+)$.

Is there a Nash equilibrium in $E_1$ for the three industries with respect to the given utility functions defined in (11–13)? If so, is it possible that the outcomes in $E_2$ transformed by the operator $A$ from this Nash equilibrium in $E_1$ constitute a Nash equilibrium in the economy $E_2$? This question can be interpreted as the existence of solutions to the split Nash equilibrium problem $SNE(G(N, S_N, f), G(M, S_M, g), A)$.

Applying optimization theory in calculus, we can show that

$$f_a(x, 2, 4) \leq f_a(1, 2, 4), \text{ for all } x \in [0, \infty);$$

$$f_b(1, y, 4) \leq f_b(1, 2, 4), \text{ for all } y \in [0, \infty);$$

$$f_c(1, 2, z) \leq f_c(1, 2, 4), \text{ for all } z \in [0, \infty).$$

It implies that the profile $(1, 2, 4) \in S_N = R_+^3$ is a Nash equilibrium for the noncooperative strategic game $G(N, S_N, f)$.

By the positive linear operator $A$, this Nash equilibrium $(1, 2, 4)$ for $G(N, S_N, f)$ can be transformed to a profile $(9, 12)$ in the game $G(M, S_M, g)$ as follows

$$A(1, 2, 4) = \begin{pmatrix} 1 & 2 & 1 \\ 2 & 1 & 2 \end{pmatrix}(1, 2, 4) = \begin{pmatrix} 9 \\ 12 \end{pmatrix}.$$

That is,

$$A_d(1, 2, 4)^T = (1\ 2\ 1)(1, 2, 4)^T = 9,$$

and
$$A_e(1, 2, 4)^T = (2\ 1\ 2)(1, 2, 4)^T = 12.$$

We can similarly prove that

$$g_d(s, 12) \leq g_d(9, 12), \text{ for all } s \in [0, \infty);$$

$$g_e(9, t) \leq g_e(9, 12), \text{ for all } t \in [0, \infty).$$

It follows that the profile $(9, 12) = (A(1, 2, 4))^T$ is a Nash equilibrium for the noncooperative strategic game $G(M, S_M, g)$. Hence $(1, 2, 4)$ is a solution to the split Nash equilibrium problem $SNE(G(N, S_N, f), G(M, S_M, g), A)$.

## 5. Split Nash equilibrium for repeating (twice) noncooperative strategic games

In game theory and in microeconomic theory, it is an important topic to study the behaviors of the equilibria when some games are repeated played. In particular, when a game is repeated, if the players select some new strategies which are created by linear transformations of the strategies in the previous action. It can be treated as a special case of split Nash equilibrium problems for two related games.

In a split Nash equilibrium problem $SNE(G(N, S_N, f), G(M, S_M, g), A)$, if the game $G(N, S_N, f)$ is exactly the same as $G(M, S_M, g)$, then we say that the game $G(N, S_N, f)$ is repeated in this problem. In this case, the linear continuous operator $A: S_N \to S_M$ can be considered a modification of the strategies when this game is repeated. The split equilibrium problem $SNE(G(N, S_N, f), G(N, S_N, f), A)$ is called a split Nash equilibrium problem for repeated game, simply denoted by $SNE(G(N, S_N, f)^2, A)$.

Let $G(N, S_N, f)$ be a $n$-person non-cooperative strategic game with concave utility function. A linear continuous operator $A: B_N \to B_N$ is said to be *convexity direction preserved* with respect to $f$ on $S_N$ if, for any given points $(u, Au), (v, Av) \in S_N \times S_N$ and for their arbitrary convex combination $w = \lambda u + (1-\lambda)v$, where $0 \leq \lambda \leq 1$, we have either that

$$F(u, w) \leq^n f(w) \text{ and } F(Au, Aw) \leq^n f(Aw) \text{ both hold,}$$

or that

$$F(v, w) \leq^n f(w) \text{ and } F(Av, Aw) \leq^n f(Aw) \text{ both hold.}$$

From Theorem 3. 3, we have

**Corollary 5.1**. *Let $G(N, S_N, f)$ be a n-person non-cooperative strategic game. Let $A: B_N \to B_N$ be a linear continuous operator. Suppose that $f_i$'s and $A$ satisfy the following conditions*:

(a1). $f_i$'s all are continuous and concave;
(a2). $AS_N = S_N$;
(a3). $A$ is convexity vector direction preserved with respect to $f$ on $S_N$.

*If there is $(t, At) \in S_N \times S_N$ such that*

$$\{(z, Az) \in S_N \times S_N: F(t, z) \leq^n f(z) \text{ and } F(At, Az) \leq^n f(Az)\} \text{ is compact,}$$

*then* $\text{SNE}(G(N, S_N, f)^2, A)$ *has a solution.*

Next we consider split Nash equilibrium problems for repeated game $\text{SNE}(G(N, S_N, f)^2, A)$, in which $A: B_N \to B_N$ is defined by a transition matrix, such as,

$$A = \begin{pmatrix} \alpha_{11} & \alpha_{12} & \cdots & \alpha_{1n} \\ \alpha_{21} & \alpha_{22} & \cdots & \alpha_{2n} \\ \cdots & \cdots & \cdots & \cdots \\ \alpha_{n1} & \alpha_{n2} & \cdots & \alpha_{nn} \end{pmatrix}, \quad (16)$$

where $\alpha_{ij} \geq 0$, for all $i, j = 1, 2, \ldots, n$ and $\sum_{j=1}^{n} \alpha_{ij} = 1$, for all $i = 1, 2, \ldots, n$.

In a split Nash equilibrium problem for repeated game $\text{SNE}(G(N, S_N, f)^2, A)$, if the linear operator $A$ is defined by a transition matrix such as (16), then it is called a *split Nash equilibrium problem of Markov type*, or simply called a *Markov split Nash equilibrium problem*. As an example of Markov split Nash equilibrium problem, we study the extended Bertrant model of price competition.

### 6. Applications to the dual extended Bertrant duopoly model of price competition

The Bertrant duopoly model of price competition (denoted by BM) is a model of oligopolistic competition that deals with two profit-maximizing firms, named by 1 and 2, in a market (see [9]). In this model, it is assumed that the two firms have constant returns to scale technologies with the same cost $c$, per unit produced. It means that the products by the two firms have the same quality.

In this section, we generalize the Bertrant duopoly model to the extended Bertrant model of price competition (denoted by EBM). Suppose that the qualities of the productions produced by the two firms possibly be different. The two firms have constant returns to scale technologies with costs $c_1 > 0$ and $c_2 > 0$, per unit produced, respectively. In the extended Bertrant Model, $c_1$ may be different from $c_2$, that explains that the qualities of the productions by these two firms may be different. In the extended Bertrant model, without loss of the generality, we assume

$$c_1 \leq c_2, \quad (17)$$

The inequality (17) means that the quality of the products in firm 1 may not be as good as that produced in firm 2.

Let $p_j$ be the price of the products by firm $j$, for $j = 1, 2$. Let $\delta(p_1, p_2)$ be the demand function in this duopoly market. Let $\delta_j(p_1, p_2)$ be the sale function for firm $j$, for $j = 1, 2$. $f$ and $\delta_j$ are assumed to be continuous functions of two variables and strictly decreasing with respect to every given variable. Suppose that there are positive numbers $\bar{p}_j$, for $j = 1, 2$, such that, for all $p_k$

$$\delta_j(p_j, p_k) \geq 0, \text{ for all } p_j \in [0, \bar{p}_j) \text{ and } \delta_j(p_j, p_k) = 0, \text{ for all } p_j \geq \bar{p}_j$$

Suppose that the socially optimal (competitive) output level in this market is strictly positive and finite for every firm
$$0 < \delta(c_1, c_2) < \infty.$$

For given prices $p_1$, $p_2$, for firms 1 and 2, respectively, the market is assumed to be clear. That is,

$$\delta(p_1, p_2) = \delta_1(p_1, p_2) + \delta_2(p_1, p_2).$$

Let $\lambda = c_1/c_2$, that defines the ratio of the qualities of the products by firm 1 to firm 2. From the assumption (17), we have $\lambda \in (0, 1]$. Considered as a noncooperative strategic game, the competition takes place as follows: The two firms simultaneously name their prices $p_1, p_2$, respectively. The sales $\delta_1(p_1, p_2)$ and $\delta_2(p_1, p_2)$ are then satisfied

$$\frac{\delta_1(p_1, p_2)}{\delta(p_1, p_2)} = \begin{cases} 0, & \text{if } p_1 > \lambda p_2 \\ \frac{c_1}{c_1 + c_2}, & \text{if } p_1 = \lambda p_2 \\ 1, & \text{if } p_1 < \lambda p_2, \end{cases} \quad (18)$$

and

$$\frac{\delta_2(p_1, p_2)}{\delta(p_1, p_2)} = \begin{cases} 1, & \text{if } p_1 > \lambda p_2 \\ \frac{c_2}{c_1 + c_2}, & \text{if } p_1 = \lambda p_2 \\ 0, & \text{if } p_1 < \lambda p_2. \end{cases} \quad (19)$$

We assume that the firms produce to order and so they incur production costs only for an output level equal to their actual sales. Therefore, for given prices $p_1, p_2$, the firm $j$ has profits

$$u_j(p_1, p_2) = (p_j - c_j)\delta_j(p_1, p_2), \text{ for } j = 1, 2. \quad (20)$$

**Theorem 6.1**. *In the extended Bertrant duopoly Model, there is a unique Nash equilibrium $(\hat{p}_1, \hat{p}_2) = (c_1, c_2)$. In this equilibrium, both firms set their prices equal to their costs, respectively*: $\hat{p}_1 = c_1$, $\hat{p}_2 = c_2$.

*Proof.* In the proof of this theorem, we always assume that the firms realistically name their prices $p_j \in [0, \bar{p}_j)$, for $j = 1, 2$. It implies $\delta(p_1, p_2) > 0$, for all prices $p_1, p_2$ realistically set by firms 1 and 2, respectively.

At first we show that $(c_1, c_2)$ is a Nash equilibrium in the extended Bertrant duopoly Model. From (20), if both firms set their prices equal to their costs, respectively, then both firms earn zero profits, that is,

$$u_j(c_1, c_2) = (c_j - c_j)\delta_j(c_1, c_2) = 0, \text{ for } j = 1, 2. \tag{21}$$

As firm 2 sets its price $c_2$, for any price $p_1$ named by firm 1, it earns profits

$$u_1(p_1, c_2) = (p_1 - c_1)f_1(p_1, c_2). \tag{22}$$

Noting $\lambda c_2 = c_1$, from (18), we have

$$\frac{\delta_1(p_1, c_2)}{\delta(p_1, c_2)} = \begin{cases} 0, & \text{if } p_1 > c_1 \\ \dfrac{c_1}{c_1 + c_2}, & \text{if } p_1 = c_1 \\ 1, & \text{if } p_1 < c_1. \end{cases}$$

By (22), it implies

$$u_1(p_1, c_2) = \begin{cases} 0, & \text{if } p_1 \geq c_1 \\ < 0, & \text{if } p_1 < c_1. \end{cases}$$

Hence

$$u_1(p_1, c_2) \leq u_1(c_1, c_2), \text{ for all price } p_1 \text{ named by firm 1}. \tag{23}$$

As firm 1 sets its price $c_1$, for any price $p_2$ named by firm 2, it earns profits

$$u_2(c_1, p_2) = (p_2 - c_2)\delta_2(c_1, p_2). \tag{24}$$

Since $\dfrac{1}{\lambda} c_1 = c_2$, from (19), we have

$$\frac{\delta_2(c_1, p_2)}{\delta(c_1, p_2)} = \begin{cases} 1, & \text{if } p_2 < c_2 \\ \dfrac{c_2}{c_1 + c_2}, & \text{if } p_2 = c_2 \\ 0, & \text{if } p_2 > c_2. \end{cases}$$

From (24), it implies

$$u_2(c_1, p_2) = \begin{cases} <0, & \text{if } p_2 < c_2 \\ 0, & \text{if } p_2 \geq c_2. \end{cases}$$

Hence

$$u_2(c_1, p_2) \leq u_2(c_1, c_2), \text{ for all price } p_2 \text{ named by firm 2.} \tag{25}$$

Combining (23) and (25), it concludes that $(c_1, c_2)$ is a Nash equilibrium in the extended Bertrant duopoly Model.

Next we prove the uniqueness of the Nash equilibrium. Assume that $(d_1, d_2)$ is also a Nash equilibrium in the extended Bertrant duopoly Model. At the given prices $(d_1, d_2)$, the earns of the two firms respectively are

$$u_1(d_1, d_2) = (d_1 - c_1)\delta_1(d_1, d_2), \tag{26}$$

$$u_2(d_1, d_2) = (d_2 - c_2)\delta_2(d_1, d_2). \tag{27}$$

Where the sale functions $\delta_1(d_1, d_2)$ and $\delta_2(d_1, d_2)$ satisfy the following equations

$$\frac{\delta_1(d_1, d_2)}{\delta(d_1, d_2)} = \begin{cases} 0, & \text{if } d_1 > \lambda d_2 \\ \dfrac{c_1}{c_1 + c_2}, & \text{if } d_1 = \lambda d_2 \\ 1, & \text{if } d_1 < \lambda d_2, \end{cases} \tag{28}$$

$$\frac{\delta_2(d_1, d_2)}{\delta(d_1, d_2)} = \begin{cases} 1, & \text{if } d_1 > \lambda d_2 \\ \dfrac{c_2}{c_1 + c_2}, & \text{if } d_1 = \lambda d_2 \\ 0, & \text{if } d_1 < \lambda d_2. \end{cases} \tag{29}$$

At first, we show $d_1 = c_1$. The proof is divided into two cases.

Case 1: Assume $d_1 < c_1$. Case 1 is also divided into two subcases with respect to $d_1$ and $\lambda d_2$.

Subcase 1.1: $d_1 > \lambda d_2$. In this case, we must have $d_2 < c_2$. Otherwise, $c_1 > d_1 > \lambda d_2 \geq \lambda c_2 = c_1$. A contradiction. Then from (27) and (29), we have

$$u_2(d_1, d_2) = (d_2 - c_2)\delta_2(d_1, d_2) < 0, \text{ if } d_1 > \lambda d_2.$$

For any price $p_2$ set by firm 2 satisfying $d_1 < \lambda p_2$, from (29), we have $\delta_2(d_1, p_2) = 0$. It implies

$$u_2(d_1, p_2) = (p_2 - c_2)\delta_2(d_1, p_2) = 0, \text{ if } d_1 < \lambda p_2.$$

It implies

$$u_2(d_1, d_2) < u_2(d_1, p_2), \text{ for all } p_2 \text{ with } d_1 < \lambda p_2. \tag{30}$$

Subcase 1.2: $d_1 \leq \lambda d_2$. In this subcase, $\delta_1(d_1, d_2) > 0$. From $d_1 < c_1$, it implies

$$u_1(d_1, d_2) < 0.$$

For any price $p_1$ set by firm 1 satisfying $p_1 > \lambda d_2$, from (28), we have $\delta_1(p_1, d_2) = 0$. It implies

$$u_1(p_1, d_2) = (p_1 - c_1)\delta_1(p_1, d_2) = 0.$$

It follows that

$$u_1(d_1, d_2) < u_1(p_1, d_2), \text{ for all } p_1 \text{ with } p_1 > \lambda d_2. \tag{31}$$

By (30) and (31), it implies that if $d_1 < c_1$, then $(d_1, d_2)$ cannot be a Nash equilibrium.

Case 2: Assume $d_1 > c_1$. Case 2 is divided into five subcases with respect to $d_1$ and $\lambda d_2$.

Subcase 2.1: $d_1 > \lambda d_2 > c_1$. In this subcase, $\delta_1(d_1, d_2) = 0$. It implies $u_1(d_1, d_2) = 0$. For any price $p_1$ named by firm 1 satisfying $c_1 < p_1 \leq \lambda d_2$, from (28), we have $\delta_1(p_1, d_2) > 0$ and

$$u_1(p_1, d_2) = (p_1 - c_1)\delta_1(p_1, d_2) > 0.$$

It implies

$$u_1(d_1, d_2) < u_1(p_1, d_2), \text{ for all } p_1 \text{ with } c_1 < p_1 \leq \lambda d_2. \tag{32}$$

Subcase 2.2: $d_1 > c_1 > \lambda d_2$. Since $\lambda c_2 = c_1$, it implies $d_2 < c_2$. In this case, from (29), we have

$$u_2(d_1, d_2) < 0 \text{ and } u_2(d_1, c_2) = 0. \tag{33}$$

Subcase 2.3: $d_1 > \lambda d_2 = c_1$. Since $\lambda c_2 = c_1$, it implies $d_2 = c_2$. If firm 2 names a price $p_2$ satisfying $d_1 > \lambda p_2 > c_1$, then, $\lambda p_2 > c_1 = \lambda d_2 = \lambda c_2$. It implies $p_2 > c_2$. From (27) and (29), we have

$$u_2(d_1, d_2) = 0 \text{ and } u_2(d_1, p_2) > 0. \tag{34}$$

Subcase 2.4: $d_1 = \lambda d_2$. In this case, we must have $d_2 > c_2$. Otherwise, $c_1 < d_1 = \lambda d_2 \leq \lambda c_2 = c_1$. A contradiction. Then from (27) and (29), we have

$$u_2(d_1, d_2) = (d_2 - c_2) \frac{c_2}{c_1 + c_2} \delta(d_1, d_2). \tag{35}$$

We take an arbitrary small number $\epsilon > 0$. Since $d_1 > \lambda(d_2 - \epsilon)$, from (27) and (29), we have

$$u_2(d_1, d_2 - \epsilon) = (d_2 - \epsilon - c_2)\delta(d_1, d_2 - \epsilon). \tag{36}$$

Noting $0 < \frac{c_2}{c_1+c_2} < 1$, and, for fixed $d_1$, $\delta(d_1, d_2)$ is continuous and strictly decreasing with respect to $d_2$, we can choose $\epsilon$ to be positive and small enough, such that $d_2 - \epsilon > c_2$ and

$$(d_2 - c_2) \frac{c_2}{c_1+c_2} \delta(d_1, d_2) < (d_2 - \epsilon - c_2)\delta(d_1, d_2 - \epsilon).$$

By (35) and (36), it implies

$$u_2(d_1, d_2) < u_2(d_1, d_2 - \epsilon), \text{ for all } \epsilon \text{ to be small enough satisfying } d_2 - \epsilon > c_2. \quad (37)$$

Subcase 2.5: $d_1 < \lambda d_2$. Then from (27) and (29), we have $u_2(d_1, d_2) = 0$. Take $p_2$ satisfying $d_1 = \lambda p_2$. Similar to subcase 2.4, we must have $p_2 > c_2$. It implies

$$u_2(d_1, p_2) = (p_2 - c_2) \frac{c_2}{c_1+c_2} \delta(d_1, p_2) > 0.$$

We have

$$u_2(d_1, d_2) < u_2(d_1, p_2), \text{ for all } p_2 \text{ satisfying } d_1 = \lambda p_2. \quad (38)$$

Combining (32–34), (37–38), it implies that if $d_1 > c_1$, then $(d_1, d_2)$ is not a Nash equilibrium. By summarizing cases 1 and 2, we must have $d_1 = c_1$.

Now we assume $d_1 = c_1$ and then show $d_2 = c_2$. The proof is divided into case I and case II.

Case I: Assume $d_2 < c_2$. Case I is also divided into two subcases with respect to $c_1$ and $\lambda d_2$.

Subcase I.1: $c_1 \geq \lambda d_2$. Then from (27) and (29), we have $u_2(c_1, d_2) = (d_2 - c_2)\delta_2(c_1, d_2) < 0$. For any price $p_2$ set by firm 2 satisfying $c_1 < \lambda p_2$, from (29), we have $\delta_2(c_1, p_2) = 0$. It implies

$$u_2(c_1, p_2) = (p_2 - c_2)\delta_2(c_1, p_2) = 0, \text{ if } c_1 < \lambda p_2.$$

We get
$$u_2(c_1, d_2) < u_2(c_1, p_2), \text{ for all } p_2 \text{ with } c_1 < \lambda p_2. \quad (39)$$

Subcase I.2: $c_1 < \lambda d_2$. In this subcase, $\delta_1(c_1, d_2) > 0$ and $u_1(c_1, d_2) = 0$. For any price $p_1$ set by firm 1 satisfying $c_1 < p_1 < \lambda d_2$, from (28), we have $\delta_1(p_1, d_2) > 0$. It implies

$$u_1(p_1, d_2) = (p_1 - c_1)\delta_1(p_1, d_2) > 0.$$

Then
$$u_1(c_1, d_2) < u_1(p_1, d_2), \text{ for all } p_1 \text{ with } c_1 < p_1 < \lambda d_2. \quad (40)$$

By (39) and (40), if $d_2 < c_2$, then $(c_1, d_2)$ cannot be a Nash equilibrium.

Case II: Assume $d_2 > c_2$. In this case, $\lambda d_2 > c_1$ and $u_1(c_1, d_2) = 0$. For any price $p_1$ set by firm 1 satisfying $c_1 < p_1 < \lambda d_2$, from (28), we have $\delta_1(p_1, d_2) > 0$. Then

$$u_1(p_1, d_2) = (p_1 - c_1)\delta_1(p_1, d_2) > 0.$$

It implies that if $d_2 > c_2$, then $(c_1, d_2)$ cannot be a Nash equilibrium. By summarizing cases I and II, we conclude that $d_2 = c_2$. □

It is well-known that the realistic model about producing-sale competition of some firms should be a dynamic model instead of static. In this paper, we consider the dual extended Bertrant duopoly model of price competition that is a repeating (twice) game of the static extended Bertrant studied in Theorem 6.1. In the dual extended Bertrant model, a firm always consider the reaction of its competitor to its price or quantity. For example, firm 1 could try to increase its profits by increasing the price from $p_1$ to $p_1'$ (never excess $p_2$), even though decreasing its sales $\delta_1(p_1', p_2)$ Meanwhile, firm 2 could try to increase its profits by decreasing the price from $p_2$ to $p_2'$ (never lower than $p_1$) and increasing its sales $\delta_2(p_1, p_2')$. Such performance can be considered as a linear transformation from $(p_1, p_2)$ to $(p_1', p_2')$: there is a 2×2 matrix $A$

$$A = \begin{pmatrix} \alpha & 1-\beta \\ 1-\alpha & \beta \end{pmatrix}, \tag{41}$$

where $0 \leq \alpha, \beta \leq 1$ such that $(p_1', p_2') = A(p_1, p_2)$.

We consider if there is a Nash equilibrium $(\hat{p}_1, \hat{p}_2)$ in the static extended Bertrant model of price competition such that the transformed (modified) prices $A(\hat{p}_1, \hat{p}_2)$ is also a Nash equilibrium in the extended Bertrant model. As a matter of fact, it is a Markov split Nash equilibrium problem for dual extended Bertrant model considered as repeating (twice) games. It is written as SNE(EBM$^2$, $A$), that is studied in the previous section. By Theorem 6.1, $(c_1, c_2)$ is the unique Nash equilibrium in the extended Bertrant duopoly Model of price competition. As a consequence of Theorem 6.1, we have

**Theorem 6.2**. *For the Markov split Nash equilibrium problem of the dual extended Bertrant duopoly Model of price competition* SNE(EBM$^2$, $A$), *we have*

(i) *If $c_1 = c_2 = c$, then, for any linear transformation $A$ defined in (41), $(c, c)$ is the unique Markov split Nash equilibrium in* SNE(EBM$^2$, $A$);

(ii) *If $c_1 < c_2$, then there exists a unique Markov split Nash equilibrium in* SNE(EBM$^2$, $A$), *$(c_1, c_2)$, only if the linear transformation $A$ is the identity*

$$A = \begin{pmatrix} 1 & 0 \\ 0 & 1 \end{pmatrix}.$$